\documentstyle[11pt,amscd]{amsart}
\newtheorem{thm}{Theorem}

\newtheorem{prop}[thm]{Proposition}

\theoremstyle{definition}
\newtheorem{defn}{Definition}

\theoremstyle{remark}

\begin{document}

\title[INVARIANTS OF FIBRED KNOTS FROM MODULI]{ Invariants of fibred knots from moduli}
\author{Hans U. Boden}

\address{Max-Planck-Institut f\"{u}r Mathematik, Gottfried-Claren-Str.\\ 26, 53225 Bonn, Germany}
\thanks{ This work was partially supported by a Rackham Grant from the University of Michigan} 
\email {hboden@@mpim-bonn.mpg.de }

\maketitle

\setcounter{page}{1}


In this paper,  invariants of fibred knots are described in terms of unitary representation spaces 
of  fundamental groups of Seifert surfaces.
Because these representation spaces admit an interpretation  as moduli spaces of parabolic bundles \cite{MS}, 
the results of \cite{B1, BH, BY} can be used to calculate the invariants. 
The main result here is Theorem 3.
The proof of this theorem involves the theory of 
parabolic bundles in an essential way, but for the sake of clarity we have suppressed their role.
These issues will be addressed in \cite{B2}.

The knot invariants, denoted $\mu_\alpha({\bold K})$ or simply $\mu_\alpha$ when the knot is understood, are defined
for fibred knots ${\bold K}$ in a rational homology sphere {\bf M} and for $\alpha \in {\bold S \bold U}_n$ a regular value of the 
map $\widetilde{\Psi}$ 
described in Definition 1.
They are a generalization of Frohman's invariants, 
$\mu_\omega,$ which are invariants of fibred knots defined for $\omega$ a regular value of $\widetilde{\Psi}$ 
in the center of ${\bold S \bold U}_n$ \cite{F}.
These have already been extended to arbitrary knots in rational homology spheres in the cases
\begin{enumerate}
\renewcommand{\labelenumi}{(\arabic{enumi})}
\item where $\omega$ is a regular value of $\widetilde{\Psi}$ in the center of ${\bold S \bold U}_n$ in \cite{FL}, and
\item where $\omega$ is any element  in the center of ${\bold S \bold U}_n$ in \cite{FN}.
\end{enumerate}
The invariants presented in this paper form a continuous family of invariants parameterized by
conjugacy classes of regular values 
$\alpha$ of $\widetilde{\Psi}.$
The most interesting aspect of this approach is the behavior of $\mu_{\alpha}$ 
as $\alpha$ is allowed to vary in the set 
$$W  \overset{\operatorname{def}}{=} {\bold S \bold U}_n /{\operatorname{Ad}}.$$

The critical values of $\Psi$ form a union of hyperplanes, giving $W$ a natural 
chamber structure. 
Theorem 3 describes the behavior of these invariants
in terms of this chamber structure.
Part (1) of Theorem 3 states that $\mu_\alpha = \mu_\alpha'$ for $\alpha$ and $\alpha'$ in the same chamber,
part (2) compares $\mu_\alpha$ and $\mu_\beta$ for
$\alpha$ and $\beta$ in adjacent chambers,
and part (3) relates $\mu_\alpha$ and $\mu_\omega$ for
$\alpha$ in the interior and $\omega$ on the boundary of $W.$
The interest in these invariants lies in their ability to detect certain irreducible
representations of the knot group: if $\mu_\alpha \neq 0$
then there exists
an irreducible unitary representation of the fundamental group of the knot complement with holonomy 
along the longitude conjugate to
$\alpha$. 

To start, we fix some notation.
Given an arbitrary group $\pi$ and a compact Lie group {\bf  G}, let 
$$\widetilde{\bold R}(\pi,{\bold G}) = {\operatorname{Hom}} ( \pi, {\bold G})$$
denote the space of homomorphisms of $\pi$ into {\bf G}.
The group {\bf G} acts on $\widetilde{\bold R}(\pi,{\bold G})$ by conjugation (denoted by Ad)
and we get the quotient
$${\bold R}(\pi,{\bold G}) = \widetilde{\bold R}(\pi,{\bold G}) /{\operatorname{Ad}} $$
which is called the 
space of representations of $\pi$ in {\bf  G}.

To simplify notation, set  $\pi = \pi_1 F,$ the fundamental group of a closed Riemann
surface $F$ of genus $g$ and $\pi^* = \pi_1 F^*,$ the fundamental
group of $F^\ast = F \setminus D^2(p),$ where $ D^2(p)$ is a small disk centered at $p \in F.$
Then $\pi$ admits the  presentation
$$\pi = \langle a_1, b_1, \ldots, a_g, b_g \mid \prod^g_{i=1} \, [a_i,b_i]\rangle, \leqno(1)$$
and $\pi^*$ is simply the free group on the $2g$ generators $a_1, b_1, \ldots, a_g, b_g.$
This last fact gives an identification 
$$\begin{array}{rl}
\widetilde{\bold R}(\pi^*,{\bold G}) &\cong \overbrace{{\bold G} \times \cdots \times {\bold G}}^{2g},\\
\rho & \mapsto (A_1, B_1, \ldots, A_g, B_g) \end{array}$$
where $A_i= \rho(a_i)$ and $B_i = \rho(b_i).$ 
Setting $\partial = \prod^g_{i=1} \, [a_i,b_i]\in \pi^*$ we  
define maps $\widetilde{\Psi}$ and $\Psi$ on $\widetilde{\bold R}$ and ${\bold R}$ by evaluation on $\partial.$

\begin{defn} Let $\widetilde{\Psi}:\widetilde{\bold R}(\pi^*,{\bold G}) \longrightarrow {\bold G}$ be defined using the
above identification and setting 
$\widetilde{\Psi} (A_1,B_1,\ldots,A_g,B_g) = \prod^g_{i=1} [A_i,B_i].$
Since 
$\widetilde{\Psi}(g \cdot \rho) = g \widetilde{\Psi}(\rho) g^{-1},$ we can define 
$\Psi:{\bold R}(\pi^*,{\bold G})\longrightarrow {\bold G}/{\operatorname{Ad}}$
and the following square commutes:
$$\begin{CD}
\widetilde{\bold R}(\pi^*,{\bold G}) @>{\widetilde{\Psi}}>> {\bold G} \\ 
@V{\operatorname{Ad}}VV                                       @VV\operatorname{Ad}V \\
{\bold R}(\pi^*,{\bold G}) @>{\Psi}>> {\bold G}/{\operatorname{Ad}}.
\end{CD}\leqno(2) $$
\end{defn}
By definition,  $\operatorname{im}(\widetilde{\Psi}) \subset [{\bold G},{\bold G}],$ the commutator
subgroup of {\bf G}.
The surjection $\pi^* \rightarrow \pi$ induces an inclusion 
${\bold R}(\pi,{\bold G})\hookrightarrow{\bold R}(\pi^*,{\bold G})$
defined by pullback, and it is clear
that 
$\rho \in \operatorname{im}(i) \Longleftrightarrow \Psi(\rho) = I.$

For the remainder of the paper, we set ${\bold G}={\bold U}_n,$ the group of unitary 
$n \times n$ complex matrices. Recall that ${\bold S}{\bold U}_n = [{\bold U}_n, {\bold U}_n],$
and denote by ${\bold Z}_n$ the center of ${\bold S}{\bold U}_n$ and
by ${\bold P \bold U}_n$ the quotient ${\bold S}{\bold U}_n / {\bold Z}_n.$  

Suppose now that {\bf K} is a fibred knot in a homology sphere {\bf M} with spanning surface 
$F^\ast$.  Let $$\varphi \,: \, F^\ast @>{\approx}>>  F^\ast$$ 
be the monodromy of the fibration, so the knot complement is given by 
$${\bold M}_{\bold K}=F^\ast \times I / \sim,$$
where $(x,0) \sim(\varphi(x),1).$ 

In \cite{F}, Frohman defines invariants 
of fibred knots by considering the Lefschetz number of the monodromy
action on certain smooth submanifolds of ${\bold R}(\pi^*, {\bold U}_n),$
namely ${\bold R}_\omega$ (defined below), where $\omega \in {\bold Z}_n$ is a regular
value of $\Psi.$ This is the case if and only if $\omega = e^{2 \pi i k/n} I,$ where
$k$ and $n$ are relatively prime.

We now describe an extension of these invariants. 
For any $\alpha \in {\bold S \bold U}_n,$
denote by $C(\alpha)$ the orbit ${\operatorname{Ad}}\cdot \alpha$ of $\alpha$ under conjugation and
set $\widetilde{\bold R}_\alpha =  \widetilde{\Psi}^{-1} (C(\alpha)).$ 
A standard result shows
that the critical points of $\widetilde{\Psi}$ are precisely the reducible
representations (see \cite{AM} or \cite{FL}). In particular, if $\alpha$ is a
regular value of $\widetilde{\Psi},$ then
$\widetilde{\Psi}^{-1}(C(\alpha))$ is smooth and consists entirely of irreducible 
{representations.}\footnote{It is enough to assume that $\widetilde{\Psi}$ is transverse to $C(\alpha),$
but a simple calculation shows that if $\widetilde{\Psi}$ is transverse to $C(\alpha),$ then
$\alpha$ is indeed a regular value.}

Define $${\bold R}_\alpha=\widetilde{\bold R}_\alpha/{\operatorname{Ad}} = \Psi^{-1}(C(\alpha)).$$
Since the adjoint action is a free ${\bold P \bold U}_n$ 
action on the irreducible representations, it follows that 
for $\alpha$ a regular value,
${\bold R}_\alpha$ is a compact oriented manifold of dimension
$$\dim {\bold R}_\alpha = (2g-2)n^2+2+\dim C(\alpha).$$ 
The invariants $\mu_\alpha$ are defined as follows.

\begin{defn}
\begin{enumerate}
\renewcommand{\labelenumi}{(\arabic{enumi})}
\item If $P$ is a compact oriented manifold and $f \, : \, P \rightarrow P,$ then the Lefschetz number of
$f$, denoted $\Lambda(f,P),$ is the algebraic intersection number of the graph of $f$ with the diagonal $\triangle(P)$ in $P \times P.$
Further, the Lefschetz polynomial is defined to be 
$$L(f,P)(t) = \sum_{0 \leq n} (-1)^n {\operatorname{Tr}}(f^* \, : \, {\operatorname{H}}^n(P,{\Bbb Q}) \rightarrow {\operatorname{H}}^n(P,{\Bbb Q})) t^n.$$
\item Let $\mu_\alpha ({\bold K})$ denote the Lefschetz fixed point number of the map on the representation
variety induced by the monodromy of the knot, i.e.
$$\mu_\alpha ({\bold K}) = \Lambda(\varphi^*, {\bold R}_\alpha).$$
\item Let $M_\alpha({\bold K},t)$ denote the Lefschetz polynomial of the map on the cohomology of
the representation variety induced by the monodromy of the knot, i.e.
$$M_\alpha({\bold K},t) = L(\varphi^*,{\bold R}_\alpha).$$
\end{enumerate}
\end{defn}
The Lefschetz fixed point theorem implies that the Lefschetz number equals 
the Lefschetz polynomial evaluated at $t=1,$ thus $\mu_\alpha = M_\alpha(1).$ 
Of course, if $\alpha$ is conjugate to $\alpha',$ then $\mu_\alpha = \mu_{\alpha'}.$ 
In the other cases, we would like to compare the invariants $\mu_\alpha$ and $\mu_\beta.$ 
One technique for doing this is to let $\alpha_t$ be a path connecting $C(\alpha)$ and $C(\beta)$ in 
$W$ and
to consider the map $\Psi$  restricted to the preimage of this path. By choosing the path
carefully (e.g. so it is transverse to the critical values of $\Psi$), Morse theory for $\Psi,$
in the sense of Bott, constructs a cobordism
between ${\bold R}_\alpha$ to ${\bold R}_\beta.$
The invariants $\mu_\alpha$ and $\mu_\beta$ can then be studied in terms of the
change in the cohomology under this cobordism.
This is the rough idea behind our approach.

To begin, identify the critical values of $\Psi$ with a union of hyperplanes in $W$ as follows.
Notice that
$W$ is an ($n-1$)-simplex.  (Choose a positive Weyl chamber in ${\bold s \bold u}_n$).
Use the following natural (but discontinuous!) coordinates on $W.$ 
Since any $\alpha \in {\bold S \bold U}_{n}$ is conjugate to a matrix of the form
$$\exp(diag(\alpha_1, \ldots, \alpha_n)) = \left( \matrix
e^{2 \pi i \alpha_1} & & 0 \\ & \ddots \\ 0 & & e^{2 \pi i \alpha_n}\endmatrix\right)$$
where $0 \leq \alpha_1 \leq \cdots \leq \alpha_n < 1$ and $\sum_{i=1}^n \alpha_i$ is an integer,
it follows that we can use $(\alpha_1, \ldots, \alpha_n)$ to give ${\operatorname{Ad}}$-invariant
coordinates to ${\bold S \bold U}_{n}.$ The resulting coordinates on $W$
are discontinuous precisely when $\alpha_1 = 0,$ which can be seen
by considering the coordinates assigned to the
path $\alpha_t=\exp(diag(t,\alpha_2-\frac{t}{n-1},\ldots, \alpha_n-\frac{t}{n-1}))$ for $t \in (-\epsilon, \epsilon).$

To deal with these discontinuities, decompose
$W=\bigcup_{k=0}^{n-1} W_k,$
where $$W_k = \{ (\alpha_1, \ldots \alpha_n)  \in W \mid  \sum_{i=1}^n \alpha_i = k \},$$
and notice that the coordinates are continuous along each $W_k.$
The hyperplanes $\alpha_1 = 0$ lie in $\partial W_k$ and 
are called {\it bad} hyperplanes.
Of course, $W$ can be reassembled by identifying these bad hyperplanes to the hyperplanes $\alpha_n = 1,$ 
which lie in $ \partial \overline{W_{k+1}},$ in the more or less obvious manner.

Suppose now  that 
$\rho : \pi_1(F^*) \rightarrow {\bold U}_n$ is a reducible
representation. Then up to conjugacy, we have $\operatorname{im}(\rho) \subset {\bold U}_{n_1} \times {\bold U}_{n_2}.$
Because $\gamma = \widetilde{\Psi}(\rho)$ is contained in the commutator subgroup of $\operatorname{im}(\rho),$
which is just ${\bold S \bold U}_{n_1} \times {\bold S \bold U}_{n_2},$ 
$\gamma$ is conjugate to a matrix in block form 
$$\left( \matrix\gamma^1 & 0 \\ 0  & \gamma^2 \endmatrix\right)$$ 
where $\gamma^i \in {\bold S \bold U}_{n_i}$
for $i=1,2.$
Writing $\gamma = \exp(diag(\gamma_1,\ldots,\gamma_n)),$
we see that $\gamma = \widetilde{\Psi}(\rho)$ for a reducible $\rho$ if and only if there is a proper subcollection
$0 \leq \gamma_{\sigma(1)} \leq \cdots \leq \gamma_{\sigma(n_1)} < 1$
with $\sum_{j=1}^{n_1} \gamma_{\sigma(j)}$ an integer.

This shows that the collection of critical values of $\Psi$ (which are the projection under Ad of the
critical values of $\widetilde{\Psi}$) are given by a union of hyperplanes
$\bigcup_\xi H_\xi$
in $W,$ some of which are {\it good} (i.e. not bad), others of which are bad.
We call the connected components of $W \setminus \bigcup_\xi H_\xi$  {\it chambers}.

Suppose $\alpha, \beta \in W_k$ are regular values of $\Psi$ in adjacent chambers 
and are separated by a good hyperplane $H$. Choose
$\gamma \in H,$ a generic point lying on the hyperplane separating $\alpha$ and $\beta$ (genericity means that $\gamma$
lies on no other hyperplane).
Denote by $\Sigma_\gamma$ the reducible representations of ${\bold R}_\gamma.$
Then with a choice of a complex structure $J$ on $F,$
the theorem of Mehta and Seshadri \cite{MS} provides an identification of each of these representation
spaces with a corresponding moduli space of semistable
parabolic bundles. 
Roughly speaking, the parabolic structure is determined by the eigenspaces and eigenvalues of the matrices ($\alpha, \beta$ or $\gamma$).
In this way, each of the representation spaces
${\bold R}_\alpha, {\bold R}_\beta,$ and ${\bold R}_\gamma$  
inherits the structure of a normal, projective variety.
Furthermore, since $\gamma$ lies on only one hyperplane, there is only one way for a representation
in ${\bold R}_\gamma$ to reduce, implying that $\Sigma_\gamma$ is smooth and is in fact 
the product of the lower dimensional representation spaces ${\bold R}_{\gamma^1} \times {\bold R}_{\gamma^2}.$
The next theorem follows from Theorem 3.1 of \cite{BH} (cf. Theorem 5.3 of \cite{BY} in the case of a bad hyperplane).

\begin{thm}  Choose $\alpha, \beta, \gamma \in W$ as above 
and $J$ a complex structure on $F.$ 
Then the representation spaces ${\bold R}_\alpha$ and ${\bold R}_\beta$ are related by a special birational transformation
(like a {\it flip} in Mori theory), i.e. there are projective maps $\Phi_\alpha$ and $\Phi_\beta$
$$\matrix {\bold R}_{\alpha} \quad \quad \quad {\bold R}_{\beta}\\
\phi_{\alpha} \!\!\!\! \searrow \; \swarrow  \!\! \phi_{\beta}\\ {\bold R}_{\gamma}
\endmatrix$$
(depending a priori on $J$) with the properties that 
\begin{enumerate}
\renewcommand{\labelenumi}{(\arabic{enumi})}
\item along  ${\bold R}_\gamma \setminus \Sigma_\gamma$, $\phi_\alpha$ and $\phi_\beta$ are isomorphisms,
\item  along $\Sigma_\gamma$, $\phi_\alpha$ and $\phi_\beta$ are
${\Bbb C}{\Bbb P}^{a}$ and ${\Bbb C}{\Bbb P}^{b}$ bundles, respectively, and 
\item  $a + b + 1 = {\operatorname{codim}}_{\Bbb C} \Sigma_\gamma$.
\end{enumerate}
\end{thm}
We can also compare the spaces ${\bold R}_\alpha$ and ${\bold R}_\omega,$ where $\alpha$ and $\omega$ are
regular values of $\Psi$ in the interior and on the boundary of $W,$ respectively. By the previous theorem,
we can assume that $\alpha$ and $\omega$ are not separated by any hyperplanes.
The following is a restatement of Proposition 3.4 of \cite{BH}.

\begin{prop}There is a natural projection $\psi : {\bold R}_{\alpha} \longrightarrow {\bold R}_{\omega}$ which is a fiber
bundle with fiber a flag variety ${\cal F}.$
\end{prop}

Theorem 1 and Proposition 2 can now be used to describe the behavior of the invariants $\mu_\alpha$ 
as $\alpha$ is allowed to vary within $W.$

\begin{thm}  Suppose {\bf K} is a knot in a homology sphere. The invariants $\mu_\alpha({\bold K})$ 
depend on $\alpha \in W$ in the following way. 
\begin{enumerate}
\renewcommand{\labelenumi}{(\arabic{enumi})}
\item If $\alpha$ and $\alpha'$ are regular values of $\Psi$ in the interior of $W$ contained in the same chamber,
then $$\mu_\alpha =\mu_{\alpha'}.$$
\item If $\alpha, \beta$ and $\gamma$ are chosen as in Theorem 1, then
$$\mu_\alpha - \mu_\beta = (\chi({\Bbb C}{\Bbb P}^{a}) - \chi({\Bbb C}{\Bbb P}^{b})) \mu_{\gamma^1} \mu_{\gamma^2}.$$
\item If $\alpha$ and $\omega$ are chosen as in Proposition 2, then
$$\mu_\alpha = \chi({\cal F}) \mu_\omega.$$
\end{enumerate}
\end{thm}

{\it Sketch of proof. \,} To prove part (1), notice that
${\bold R}_\alpha$ is diffeomorphic to ${\bold R}_{\alpha'}.$
This follows because $\alpha$ and $\alpha'$ are connected by a path $\alpha_t$
which misses the hyperplanes, hence $\Psi^{-1}(\alpha_t)$ is a product.
The map $\varphi^*$ acts on this product and preserves each fiber, hence the family of maps
$\varphi^*_t$ on ${\bold R}_{\alpha_t}$ describes, after pullback to ${\bold R}_{\alpha},$
an isotopy from $\varphi^*_0$ to the pullback of $\varphi^*_1.$
Part (1) follows since Lefschetz numbers are invariants of the isotopy class of a map.

For part (2), notice that a direct consequence of Theorem 1 is the identity
$${\operatorname{P}}_t({\bold R}_\alpha) ={\operatorname{P}}_t({\bold R}_\beta)+ ({\operatorname{P}}_t({\Bbb C}{\Bbb P}^{a}) -{\operatorname{P}}_t({\Bbb C}{\Bbb P}^{b}))
{\operatorname{P}}_t(\Sigma_\gamma), \leqno(*)$$
where ${\operatorname{P}}_t$ denotes the Poincar\'{e} polynomials taken with 
${\Bbb Z}$ coefficients
(cf. Corollary 3.2 of \cite{BH}).
We would like to see that this formula also 
holds on the level of the Lefschetz polynomials, i.e. that
$$M_\alpha(t) - M_\beta (t) =
({\operatorname{P}}_t({\Bbb C}{\Bbb P}^{a}) - {\operatorname{P}}_t({\Bbb C}{\Bbb P}^{b})) M_{\gamma^1}(t) M_{\gamma^2}(t).\leqno(**)$$
The direct route is to give a general algorithm for the computation of
$M_\alpha(t)$ from which one can conclude ($**$).

To explain this algorithm, we describe first the Atiyah-Bott-Nitsure method to calculate
${\operatorname{P}}_t({\bold R}_\alpha).$
Consider the case $\omega \in {\bold Z}_n.$
By \cite{AB}, we see that
${\operatorname{P}}_t({\bold R}_{\omega})$ is described as a difference of two terms,
both of which are infinite series
in $t,$ the first being the cohomology of the classifying space of the gauge group
(which is determined by the rank and genus, cf. Theorem 2.15 of \cite{AB})
and the second admitting an expression as a power series whose coefficients are polynomial
functions in $(1+t)^{2g},$
the Poincar\'{e} polynomial of the Jacobian of the surface
(see formulas (11.1) -- (11.3) of \cite{AB}).
This formula results from an equivariantly perfect stratification on the space of
holomorphic structures on a fixed topological bundle over the Riemann surface $(F, J)$.
 
Frohman adapted these ideas in order to compute $\mu_\omega$ in \cite{F}.
He expressed the Lefschetz polynomial $M_\omega(t)$
as a difference of two Lefschetz traces (which are formal power series in $t$),
the first of which is determined by the rank and genus and the second of
which admits an expression as a power series in $t$ whose coefficients
are polynomial functions in $c(t)= L(F^*,\varphi^*),$
the Lefschetz polynomial of the monodromy acting on the Jacobian.
Of course, $c(t)$
is just the Alexander polynomial and
evaluating  $M_\omega(t)$ at $t=1$ gives a formula for the Lefschetz number $\mu_\omega$
in terms of the derivatives
of the Alexander polynomial of the knot evaluated at $t=1$ (see Theorems 1.6 and 3.14 of \cite{F}).

The Atiyah-Bott procedure was extended to parabolic bundles in 
\cite{N} and used to give explicit computations
of ${\operatorname{P}}_t ({\bold R}_{\alpha})$ for low rank \cite{B1}.
In general, the procedure expresses the Poincar\'{e} polynomial of
${\bold R}_{\alpha}$ as the difference of two infinite series, the first being
determined by the rank, genus, and conjugacy class of $\alpha$ (cf. formulas (13-15) of \cite{B1}),
and the second being a power series whose coefficients are 
polynomial functions in $(1+t)^{2g}.$ 
As in the non-parabolic case, this formula results from an equivariantly perfect stratification
on the space of holomorphic structures on a fixed topological parabolic bundle over $(F, J).$
 
This stratification leads to an algorithm for calculating $M_\alpha(t).$ 
As mentioned before, such a formula exists on
the level of cohomology, and
the crucial point is to see that the map induced by $\varphi$ 
preserves the stratification. This is only true after taking into account 
the effect of $\varphi$ on 
the choice of complex structure $J$ on $F.$ 

Having established such a formula, it is relatively easy to see
which unstable strata are created and destroyed as the monodromy condition is allowed to vary.
(This is another manifestation of equation ($*$).)
Consequently, equation ($**$) is seen to follow and this implies part (2).

To prove part (3),
consider the case  $\omega \in {\bold Z}_n.$
Define $\widehat{\bold R}_\omega$ to be the quotient of $\widetilde{\bold R}_\omega$
by conjugation by the maximal torus of ${\bold S}{\bold U}_n.$ 
It is elementary to see that $\widehat{\bold R}_\omega$ is diffeomorphic to
${\bold R}_\alpha$ in this case, and the proof of part (1) shows that 
$\mu_\alpha$ is equal to the Lefschetz number
of $\varphi^*$ on $\widehat{\bold R}_\omega$.
The natural projection
$\psi : \widehat{\bold R}_\omega \longrightarrow {\bold R}_\omega$  is a fiber bundle
with fiber ${\cal F}$, the flag variety of full flags in ${\Bbb C}^n,$ and part (3) now follows by
observing that the following diagram commutes
$$\CD
\widehat{\bold R}_{\omega}          @>{\varphi^*}>>               \widehat{\bold R}_{\omega} \\
@V{\psi}VV                                    @VV{\psi}V\\
{\bold R}_{\omega}        @>{\varphi^*}>>        {\bold R}_{\omega}.
\endCD$$

An alternative approach to proving part (3) is to use the identity
$${\operatorname{P}}_t({\bold R}_{\alpha}) = {\operatorname{P}}_t({\cal F}) {\operatorname{P}}_t({\bold R}_{\omega})$$
together with an argument similar to that given for part (2).

\medskip \noindent
{\it Remarks. \,}  
It seems likely that one could define invariants $\mu_\alpha$
for arbitrary (non-fibred) knots using either the approach of \cite{FL} or that of \cite{FN}.
The results of \cite{BH} are  promising for the computation
of the latter invariants. A direct consequence of Theorem 1 is that one
of the maps $\phi_\alpha, \phi_\beta$ is a small resolution,
hence the 
intersection homology of ${\bold R}_\gamma$ is given by the homology of either ${\bold R}_\alpha$ or ${\bold R}_\beta.$
More generally, if $\gamma$ lies on more than one hyperplane, 
then ${\bold R}_\gamma$
may have a very complicated singular locus $\Sigma_\gamma,$ but Conjecture 4.8 of \cite{BH} predicts the existence
of a small resolution from a nearby non-singular ${\bold R}_\alpha.$ 
This conjecture is true if, for example, the hyperplanes containing $\gamma$ are in general position.

\vspace{10mm}

\noindent{\bf Aknowledgements.} 
I would like to thank Charles Frohman, Arthur Greenspoon,
Andrew Nicas, Richard Wentworth, and K\^oji Yokogawa
for their advice and correspondence.

\bibliographystyle{amsplain}

\end{document}